\documentclass[a4paper,pdftex,reqno,10pt,english]{amsart}

\usepackage{geometry}
\usepackage{graphicx}
\usepackage{courier}
\usepackage[usenames]{color}
\usepackage{hyperref}
\usepackage{amssymb,amsmath,amsthm}

\newcommand{\Tsm}{\hspace*{0.31cm}}
\newcommand{\n}{\!}
\newtheorem{algo}{Algorithm}
\newtheorem{theorem}{Theorem}

\newtheorem{lemma}{Lemma}

\begin{document}

\title{Enumeration of integral tetrahedra}

\author{Sascha Kurz}

\address{Sascha Kurz\\Department of Mathematics, Physic and Informatics\\University of Bayreuth\\Germany}
\email{sascha.kurz@uni-bayreuth.de}

\begin{abstract}
 We determine the numbers of integral tetrahedra with diameter $d$ up to isomorphism
 for all $d\le 1000$ via computer enumeration. Therefore we give an algorithm that enumerates the
 integral tetrahedra with diameter at most $d$ in $O(d^5)$ time and an algorithm that can check the
 canonicity of a given integral tetrahedron with at most $6$ integer comparisons. For the number
 of isomorphism classes of integral $4\times 4$ matrices with diameter $d$
 fulfilling the triangle inequalities we derive an exact formula.
\end{abstract}


\maketitle

\section{Introduction}

Geometrical objects with integral side lengths have fascinated mathematicians for ages. A very simple geometric object is an $m$-dimensional simplex. Recently an intriguing bijection between $m$-dimensional simplices with edge lengths in $\{1,2\}$ and the partitions of $m+1$ was discovered \cite{bijection}. So far, for $m$-dimensional simplices with edge lengths in $\{1,2,3\}$ no formula is known and exact numbers are obtained only up to $m=13$ \cite{phd_kurz}. Let us more generally denote by $\alpha(m,d)$ the number of non-isomorphic $m$-dimensional simplices with edge lengths in $\{1,\dots,d\}$ where at least one edge has length $d$. We also call $d$ the diameter of the simplex. The known results, see i.e. \cite{phd_kurz}, are, besides some exact numbers,
\begin{eqnarray*}
  \alpha(1,d) &=& 1,\\
  \alpha(2,d) &=& \left\lfloor\frac{d+1}{2}\right\rfloor\left\lfloor\frac{d+2}{2}\right\rfloor=
                  \left\lfloor\frac{(d+1)^2}{4}\right\rfloor,\,\quad\quad\quad\quad\quad\quad[A002620]\\
  \alpha(m,1) &=& 1,\\
  \alpha(m,2) &=& p(m+1)-1,\;\quad\quad\quad\quad\quad\quad\quad\quad\quad\quad\quad\quad\quad\quad
  [A000065]
\end{eqnarray*}
where $p(m+1)$ denotes the number of partitions [A000041] of $m+1$. The aim of this article is the determination of the number of non-isomorphic integral tetrahedra $\alpha(3,d)$.

Besides an intrinsic interest in integral simplices their study is useful in field of integral point sets. These are sets of $n$ points in the $m$-dimensional Euclidean space $\mathbb{E}^m$ with pairwise integral distances. Applications for this combinatorial structure involving geometry and number theory are imaginable in radio astronomy (wave lengths), chemistry (molecules), physics (energy quantums), robotics, architecture, and other fields, see \cite{integral_distances_in_point_sets} for an overview. We define the largest occurring distance of an integral point set $\mathcal{P}$ as its diameter. From the combinatorial point of view there is a natural interest in the determination of the minimum possible diameter $d(m,n)$ for given parameters $m$ and $n$ \cite{integral_distances_in_point_sets,minimum_diameter,hab_kemnitz,paper_carpets,phd_kurz,characteristic,paper_laue,kurz_wassermann,dipl_piepmeyer,note_on_integral_distances}. In most cases exact values of $d(m,n)$ are obtained by an exhaustive enumeration of integral point sets with diameter $d\le d(m,n)$. A necessary first step for the enumeration of $m$-dimensional integral point sets is the enumeration of $m$-dimensional integral simplices. Hence there is a need for an efficient enumeration algorithm.

Another application of integral tetrahedra concerns geometric probabilities. Suppose you are given a symmetric $3\times 3$ matrix $\Delta_2$ with entries being equi-distributed in $[0,1]$ and zeros on the main diagonal. The probability $\mathcal{P}_2$ that $\Delta_2$ is the distance matrix of a triangle in the Euclidean metric can be easily determined to be $\mathcal{P}_2=\frac{1}{2}$. As a generalization we ask for the probability $\mathcal{P}_m$ of a similar defined $(m+1)\times(m+1)$ matrix $\Delta_m$ being the distance matrix of an $m$-dimensional simplex in the Euclidean metric. To analyze the question for $m=3$ we consider a discretization and obtain $\mathcal{P}_3=\lim\limits_{d\to\infty}\frac{4\cdot\alpha(3,d)}{d^5}$.

\noindent
Our main results are the determination of $\alpha(3,d)$ for $d\le 1000$,

\begin{theorem}
  \label{thm_alpha_hat}
   The number $\hat{\alpha}_\le(d,3)$ of  symmetric $4\times 4$ matrices with entries in
  $\{1,\dots,d\}$ fulfilling the triangle inequalities is given by

  \vspace*{-5mm}

  \begin{eqnarray*}
    \hat{\alpha}_{\le}(d,3)=
    \left\lbrace
      \begin{array}{rl}
        \frac{17d^6+425d^4+1628d^2}{2880} & \mbox{for } d\equiv 0\mod 2; \\
        \frac{17d^6+425d^4+1763d^2+675}{2880} & \mbox{for } d\equiv 1\mod 2.
      \end{array}
    \right.
  \end{eqnarray*}
  If we additionally request a diameter of exactly $d$ we have

  \vspace*{-5mm}

  \begin{eqnarray*}
    \hat{\alpha}(d,3)=
    \left\lbrace
      \begin{array}{rl}
         \frac{34d^5-85d^4+680d^3-962d^2+1776d-960}{960} & \mbox{for } d\equiv 0\mod 2; \\ 
         \frac{34d^5-85d^4+680d^3-908d^2+1722d-483}{960} & \mbox{for } d\equiv 1\mod 2,
      \end{array}
    \right.
  \end{eqnarray*}
\end{theorem}

\begin{theorem}
  \label{thm_bounds}
  $$
    0.090\le\mathcal{P}_3\le 0.111,
  $$
\end{theorem}

\noindent
and the enumeration algorithms of Section \ref{sec_orderly} and Section \ref{sec_can}, which allows us to
enumerate integral tetrahedra with diameter at most $d$ in time $O(d^5)$ and to check a $4\times 4$-matrix for canonicity using at most $6$ integer comparisons.

\section{Number of integral tetrahedra} 

Because a symmetric $4\times 4$-matrix with zeros on the diagonal has six independent non-zero values there are $d^6$ labeled integral such matrices with diameter at most $d$. To obtain the number $\overline{\alpha}_\le(d,3)$ of unlabeled matrices we need to apply the following well known Lemma:

\begin{lemma}{(Cauchy-Frobenius, weighted form)}\\
  \label{lemma_cauchy_frobenius}
  \noindent
  Given a group action of a finite group $G$ on a set $S$ 
  and a map $w:S\longrightarrow R$ from $S$ into a commutative ring $R$ containing 
  $\mathbb{Q}$ as a subring. If $w$ is constant on the orbits of $G$ on $S$, then we have,
  for any transversal $\mathcal{T}$ of the orbits:
  $$
    \sum_{t\in\mathcal{T}}w(t)=\frac{1}{|G|}\sum_{g\in G}\sum_{s\in S_g}w(s)
  $$
  where $S_g$ denotes the elements of $S$ being fixed by $g$, i.e.
  $$
    S_g=\lbrace s\in S|s=gs\rbrace \,.
  $$
\end{lemma}

For a proof, notation and some background we refer to \cite{0951.05001}. Applying the lemma yields:
\begin{lemma}
  \label{lemma_tuples}
  $$\overline{\alpha}_\le(d,3)=\frac{d^6+9d^4+14d^2}{24}$$
  and
  $$
    \overline{\alpha}(d,3)=\overline{\alpha}_\le(d,3)-\overline{\alpha}_\le(d-1,3)=
    \frac{6d^5-15 d^4+56d^3-69d^2+70d-24}{24}.
  $$
\end{lemma}

As geometry is involved in our problem we have to take into account some properties of Euclidean spaces. In the Euclidean plane $\mathbb{E}^2$ the possible occurring triples of distances of triangles are completely characterized by the triangle inequalities. In general there is a set of inequalities using the so called Cayley-Menger determinant to characterize whether a given symmetric $(m+1)\times(m+1)$ matrix with zeros on the diagonal is a distance matrix of an $m$-dimensional simplex \cite{54.0622.02}. For a tetrahedron with distances $\delta_{i,j}$, $0\le i\le j<4$, the inequality
\begin{equation}
  \label{eq_tetrahedron}
  CMD_3=\left|
   \begin{array}{ccccc}
    0 & \delta_{0,1}^2 & \delta_{0,2}^2 & \delta_{0,3}^2 & 1 \\ 
    \delta_{1,0}^2 & 0 & \delta_{1,2}^2 & \delta_{1,3}^2 & 1 \\ 
    \delta_{2,0}^2 & \delta_{2,1}^2 & 0 & \delta_{2,3}^2 & 1 \\ 
    \delta_{3,0}^2 & \delta_{3,1}^2 & \delta_{3,2}^2 & 0 & 1 \\ 
    1 & 1 & 1 & 1 & 0
  \end{array} 
  \right|
   >0
\end{equation}
has to be fulfilled besides the triangle inequalities.

In a first step we exclusively consider the triangle inequalities for $m=3$ and count the number $\hat{\alpha}_\le(d,3)$ of non-isomorphic symmetric $4\times 4$ matrices with entries in $\{1,\dots,d\}$ fulfilling the triangle inequalities.

\noindent
{\it Proof of Theorem \ref{thm_alpha_hat}.}\\
Counting labeled symmetric $4\times 4$ matrices with entries in $\{1,\dots,d\}$ fulfilling the triangle inequalities is equivalent to determining integral points in a six-dimensional polyhedron. Prescribing the complete automorphism group results in some further equalities and an application of the inclusion-exclusion principle. Thus, after a lengthy but rather easy computation we can apply Lemma \ref{lemma_cauchy_frobenius} and obtain
$$
  24\hat{\alpha}_{\le}(d,3)=3\cdot\left\lceil\frac{4d^4+5d^2}{12}\right\rceil
  +6\cdot\frac{37d^4-18d^3+20d^2-21d+(36d^2+42)
  \left\lceil\frac{d}{2}\right\rceil}{96}
$$
$$
  +\left\lceil\frac{34d^6+55d^4+136d^2}{240}\right\rceil+6\cdot\left(d^2-d\left\lceil\frac{d}{2}
  \right\rceil+\left\lceil\frac{d}{2}\right\rceil^2\right)
  +8\cdot\left(d^2-d\left\lceil\frac{d}{2}\right\rceil+\left\lceil\frac{d}{2}\right\rceil^2\right),
$$
which can be modified to the stated formulas. \hfill{$\square$}

\begin{table}[!ht]
  \begin{center}
    \begin{tabular}{||rr|rr|rr|rr|rr||}
      \hline
      \hline
      \n$d$\n&\n$\alpha(d,3)$\n&\n$d$\n&\n$\alpha(d,3)$\n&\n$d$\n&\n$\alpha(d,3)$\n&\n$d$\n&
      \n$\alpha(d,3)$\n&\n$d$\n&
      \n$\alpha(d,3)$\n\\
      \hline  
      \n 1\n&\n     1\n&\n26\n&\n 305861\n&\n51\n&\n 8854161\n&\n 76\n&\n 65098817\n&\n120\n&\n
      639349793\n\\ 
      \n 2\n&\n     4\n&\n27\n&\n 369247\n&\n52\n&\n 9756921\n&\n 77\n&\n 69497725\n&\n140\n&\n
      1382200653\n\\
      \n 3\n&\n    16\n&\n28\n&\n 442695\n&\n53\n&\n10732329\n&\n 78\n&\n 74130849\n&\n160\n&\n
      2695280888\n\\
      \n 4\n&\n    45\n&\n29\n&\n 527417\n&\n54\n&\n11783530\n&\n 79\n&\n 79008179\n&\n180\n&\n
      4857645442\n\\
      \n 5\n&\n   116\n&\n30\n&\n 624483\n&\n55\n&\n12916059\n&\n 80\n&\n 84138170\n&\n200\n&\n
      8227353208\n\\
      \n 6\n&\n   254\n&\n31\n&\n 735777\n&\n56\n&\n14133630\n&\n 81\n&\n 89532591\n&\n220\n&\n
      13251404399\n\\
      \n 7\n&\n   516\n&\n32\n&\n 861885\n&\n57\n&\n15442004\n&\n 82\n&\n 95198909\n&\n240\n&\n
      20475584436\n\\
      \n 8\n&\n   956\n&\n33\n&\n1005214\n&\n58\n&\n16845331\n&\n 83\n&\n101149823\n&\n260\n&\n
      30554402290\n\\
      \n 9\n&\n  1669\n&\n34\n&\n1166797\n&\n59\n&\n18349153\n&\n 84\n&\n107392867\n&\n280\n&\n
      44260846692\n\\
      \n10\n&\n  2760\n&\n35\n&\n1348609\n&\n60\n&\n19957007\n&\n 85\n&\n113942655\n&\n300\n&\n
      62496428392\n\\
      \n11\n&\n  4379\n&\n36\n&\n1552398\n&\n61\n&\n21678067\n&\n 86\n&\n120807154\n&\n320\n&\n
      86300970558\n\\
      \n12\n&\n  6676\n&\n37\n&\n1780198\n&\n62\n&\n23514174\n&\n 87\n&\n127997826\n&\n340\n&\n
      116862463817\n\\
      \n13\n&\n  9888\n&\n38\n&\n2033970\n&\n63\n&\n25473207\n&\n 88\n&\n135527578\n&\n360\n&\n
      155526991341\n\\
      \n14\n&\n 14219\n&\n39\n&\n2315942\n&\n64\n&\n27560402\n&\n 89\n&\n143409248\n&\n380\n&\n
      203808692441\n\\
      \n15\n&\n 19956\n&\n40\n&\n2628138\n&\n65\n&\n29783292\n&\n 90\n&\n151649489\n&\n400\n&\n
      263399396125\n\\
      \n16\n&\n 27421\n&\n41\n&\n2973433\n&\n66\n&\n32145746\n&\n 91\n&\n160268457\n&\n420\n&\n
      336178761892\n\\
      \n17\n&\n 37062\n&\n42\n&\n3353922\n&\n67\n&\n34657375\n&\n 92\n&\n169272471\n&\n440\n&\n
      424224122232\n\\
      \n18\n&\n 49143\n&\n43\n&\n3773027\n&\n68\n&\n37322859\n&\n 93\n&\n178678811\n&\n460\n&\n
      529820175414\n\\
      \n19\n&\n 64272\n&\n44\n&\n4232254\n&\n69\n&\n40149983\n&\n 94\n&\n188496776\n&\n480\n&\n
      655468974700\n\\
      \n20\n&\n 82888\n&\n45\n&\n4735254\n&\n70\n&\n43145566\n&\n 95\n&\n198743717\n&\n500\n&\n
      803900006590\n\\
      \n21\n&\n105629\n&\n46\n&\n5285404\n&\n71\n&\n46318399\n&\n 96\n&\n209427375\n&\n520\n&\n
      978079728301\n\\
      \n22\n&\n133132\n&\n47\n&\n5885587\n&\n72\n&\n49673679\n&\n
      97\n&\n220570260\n&\n540\n&\n1181221582297\n\\
      \n23\n&\n166090\n&\n48\n&\n6538543\n&\n73\n&\n53222896\n&\n
      98\n&\n232180129\n&\n560\n&\n1416796092768\n\\
      \n24\n&\n205223\n&\n49\n&\n7249029\n&\n74\n&\n56969822\n&\n
      99\n&\n244275592\n&\n580\n&\n1688540496999\n\\
      \n25\n&\n251624\n&\n50\n&\n8019420\n&\n75\n&\n60926247\n&\n100\n&\n256866619\n&\n600\n&
      \n2000468396580\n\\
      \hline
      \hline
    \end{tabular}
    \caption{Number $\alpha(d,3)$ of integral tetrahedra with diameter $d$ - part 1.}
    \label{table_number_tetrahedra}
  \end{center}
\end{table}

In addition to this proof we have verified the stated formula for $d\le 500$ via a computer enumeration.
We remark that $\frac{\hat{\alpha}_{\le}(d,3)}{\overline{\alpha}_{\le}(d,3)}$ and $\frac{\hat{\alpha}(d,3)}{\overline{\alpha}(d,3)}$ 
tend to $\frac{17}{120}=0.141\overline{6}$ if $d\rightarrow\infty$. Moreover we were able to obtain an exact formula for $\hat{\alpha}(d,3)$ because the Cayley-Menger determinant 
$$
  CMD_2=\left|
   \begin{array}{cccc}
    0 & \delta_{0,1}^2 & \delta_{0,2}^2 & 1 \\ 
    \delta_{1,0}^2 & 0 & \delta_{1,2}^2 & 1 \\ 
    \delta_{2,0}^2 & \delta_{2,1}^2 & 0 & 1 \\ 
    1 & 1 & 1 & 0
  \end{array} 
  \right|
$$
for dimension $m=2$ can be written as  
$$
  CMD_2=-(\delta_{0,1}+\delta_{0,2}+\delta_{1,2})(\delta_{0,1}+\delta_{0,2}-\delta_{1,2})
  (\delta_{0,1}-\delta_{0,2}+\delta_{1,2})(-\delta_{0,1}+\delta_{0,2}+\delta_{1,2}).
$$
Thus $CMD_2<0$ is equivalent to the well known linear triangle inequalities $\delta_{0,1}+\delta_{0,2}>\delta_{1,2}$, $\delta_{0,1}+\delta_{1,2}>\delta_{0,2}$ and $\delta_{0,2}+\delta_{1,2}>\delta_{0,1}$. Unfortunately for $m\ge 3$ the Cayley-Menger
determinant is irreducible \cite{irreducible} and one cannot simplify $(-1)^{m+1}CMD_m>0$ into a set of inequalities of lower degree. So we are unable to apply the same method to derive an analytic formula for $\alpha(d,3)$.

\begin{table}[ht]
  \begin{center}
    \begin{tabular}{||rr|rr|rr||}
      \hline
      \hline
        $d$ & $\alpha(d,3)$ & $d$ & $\alpha(d,3)$ & $d$ & $\alpha(d,3)$ \\
      \hline
      620 & 2356880503873 & 760 &  6523288334629 &  900 & 15192308794063 \\
      640 & 2762373382787 & 780 &  7428031732465 &  920 & 16957109053082 \\
      660 & 3221850132593 & 800 &  8430487428682 &  940 & 18882231158104 \\
      680 & 3740530243895 & 820 &  9538364312059 &  960 & 20978358597822 \\
      700 & 4323958989350 & 840 & 10759766492473 &  980 & 23256639532080 \\
      720 & 4978017317882 & 860 & 12103204603044 & 1000 & 25728695195597 \\
      740 & 5708932993276 & 880 & 13577602128303 &      &                \\
      \hline
      \hline
    \end{tabular}
    \label{table_number_tetrahedra_2}
    \caption{Number $\alpha(d,3)$ of integral tetrahedra with diameter $d$ - part 2.}
  \end{center}
\end{table}

\begin{lemma}
  We have $\alpha(3,d)\in \Omega(d^5)$, $\alpha(3,d)\in O(d^5)$, $\alpha_{\le}(3,d)\in \Omega(d^6)$, and
  $\alpha_{\le}(3,d)\in O(d^6)$.
\end{lemma}
\begin{proof}
  The upper bounds are trivial since they also hold for symmetric matrices with integer values at most
  $d$ and zeros on the diagonal. For the lower bounds we consider six-tuples
  $\delta_{0,1}\in[d,d(1-\varepsilon))$, $\delta_{0,2}\in[d(1-\varepsilon),d(1-2\varepsilon))$,
  $\delta_{1,2}\in[d(1-2\varepsilon),d(1-3\varepsilon))$,
  $\delta_{0,3}\in[d(1-3\varepsilon),d(1-4\varepsilon))$,
  $\delta_{1,3}\in[d(1-4\varepsilon),d(1-5\varepsilon))$, and
  $\delta_{2,3}\in[d(1-5\varepsilon),d(1-6\varepsilon))$.
  For each $\varepsilon$ there are $\Omega(d^6)$ non-isomorphic matrices. If $\varepsilon$ is suitable
  small then all these matrices fulfill the triangle conditions and inequality \ref{eq_tetrahedron}.
\end{proof}

In general we have $\alpha_{\le}(m,d)\in\Omega(d^{m(m+1)/2})$, $\alpha_{\le}(m,d)\in O(d^{m(m+1)/2})$,
$\alpha(m,d)\in\Omega(d^{m(m+1)/2-1})$, and $\alpha(m,d)\in O(d^{m(m+1)/2-1})$.

In Section \ref{sec_orderly} and Section \ref{sec_can} we give an algorithm to obtain $\alpha(d,3)$ via implicit computer enumeration. Some of these computed values are given in Table \ref{table_number_tetrahedra} and Table  \ref{table_number_tetrahedra_2}. For a complete list of $\alpha(d,3)$ for $d\le 1000$ we refer to \cite{own_hp}. This amounts to
$$
  \alpha_{\le}(1000,3)=4299974867606266\approx 4.3\cdot 10^{15}.
$$

\section{Bounds for $\mathbf{\mathcal{P}_3}$}
\label{sec_estimation}

In this section we give bounds for the probability $\mathcal{P}_3$ that $\Delta_3$ is the distance matrix of a tetrahedron in the $3$-dimensional Euclidean space $\mathbb{E}^3$, where $\Delta_3$ is a symmetric $4 \times 4$ matrix with zeros on the main diagonal and the remaining entries being equi-distributed in $[0,1]$. Therefore we consider a discretization. Let $d$ be a fixed number. We consider the $d^6$ six-dimensional
cubes $\mathcal{C}_{i_1,\dots,i_6}:=\times_{j=1}^6\left[\frac{i_j}{d},\frac{i_j+1}{d}\right]
\subseteq[0,1]^6$. For every cube $\mathcal{C}$ it is easy to decide whether every point of $\mathcal{C}$ fulfills the triangle conditions, no points of $\mathcal{C}$ fulfill the triangle conditions, or both cases occur. For inequality \ref{eq_tetrahedron} we have no explicit test but we are able to compute a lower bound $\underline{CMD_3}(\mathcal{C})$ and an upper bound $\overline{CMD_3}(\mathcal{C})$, so that we have
$$\underline{CMD_3}(\mathcal{C})\le CMD_3(x)\le\overline{CMD_3}(\mathcal{C})\text{ for all } x\in\mathcal{C}.$$
Thus for some cubes $\mathcal{C}$ we can decide that all $x\in\mathcal{C}$ correspond to a tetrahedron. We denote this case by $\Xi(\mathcal{C})=1$. If no $x\in\mathcal{C}$ corresponds to a tetrahedron we set $\Xi(\mathcal{C})=-1$. In all other cases we define $\Xi(\mathcal{C})=0$. With this we obtain for all $d\in\mathbb{N}$ the following bounds:

\begin{lemma}
  $$
    \sum_{\mathcal{C}\,:\,\Xi(\mathcal{C})=1}\frac{1}{d^6}
    \le\mathcal{P}_3\le 1-\sum_{\mathcal{C}\,:\,\Xi(\mathcal{C})=-1}\frac{1}{d^6}.
  $$
\end{lemma}

Thus we have a method to obtain bounds on $\mathcal{P}_3$ using computer calculations. For the actual computation we use two further speed ups. We can take advantage of symmetries and use an adaptive strategy: We start with a small value of $d$ and subdivide cubes $\mathcal{C}$ with $\Xi(\mathcal{C})=0$ recursively into $8$ smaller cubes. After a computer calculation we obtain
$$
  0.090\le\mathcal{P}_3\le 0.111,
$$
which proves Theorem \ref{thm_bounds}. Clearly Theorem \ref{thm_bounds} can be improved by simply letting the computers work for a longer time or by using a computing cluster, but the convergence of our approach seems to be rather slow. An enhanced check whether a cube $C$ fulfills inequality (\ref{eq_tetrahedron}) would be very useful.

Good estimates for $\mathcal{P}_3$ can be obtained by considering the values $\alpha(3,d)$ in the following way. At first we consider the probability $\tilde{\mathcal{P}}_3$ being defined as $\mathcal{P}_3$ where additionally $\delta_{0,1}=1$.

\begin{lemma}
  $$\tilde{\mathcal{P}}_3=\mathcal{P}_3.$$
\end{lemma}
\begin{proof}
  The problem of determining $\mathcal{P}_3$ or $\tilde{\mathcal{P}}_3$ is an integration problem.
  Due to symmetry we only need to consider the domain where $\delta_{0,1}$ is the maximum. For every
  $\delta_{0,1}\in(0,1]$ there is a probability $p(\delta_{0,1})$ that $\delta_{0,1},\dots,\delta_{2,3}$
  are distances of a tetrahedron where $\delta_{0,2},\dots,\delta_{2,3}\in(0,\delta_{0,1}]$ are
  equi-distributed random variables. Since $p(\delta_{0,1})$ is constant we can conclude the stated
  equation.
\end{proof}

\begin{lemma}
$$
  \mathcal{P}_3=\lim_{d\to\infty}\frac{4\cdot\alpha(d,3)}{d^5}.
$$
\end{lemma}
\begin{proof}
  We consider a modified version of the algorithm described above to obtain exact bounds on
  $\tilde{\mathcal{P}}_3$. As already mentioned, the triangle inequalities alone define a
  five-dimensional polyhedron. Since determinants are continuous $CMD_3=0$ defines a
  smooth surface and so the volume of all cubes $\mathcal{C}$ with $\Xi(\mathcal{C})=0$ converges
  to zero. Thus substituting $\Xi(\mathcal{C})$ by the evaluation of $\Xi$ in an arbitrary corner
  of $\mathcal{C}$ yields the correct value for $\tilde{\mathcal{P}}_3=\mathcal{P}_3$ for
  $d\to\infty$. Since there are at most $O(d^4)$ six-tuples $(d,i_2,\dots,i_6)$, $i_j\in\{1,\dots,d\}$
  with non-trivial automorphism group we can factor out symmetry and conclude the stated result.
\end{proof}

Using Lemma \ref{lemma_tuples} and Theorem \ref{thm_alpha_hat} we can modify this to
$$
  \mathcal{P}_3=\lim_{d\to\infty}\frac{\alpha(d,3)}{\overline{\alpha}(d,3)}\le\lim_{d\to\infty}
  \frac{\hat{\alpha}(d,3)}{\overline{\alpha}(d,3)}=\frac{17}{120}=0.141\overline{6}.
$$
Heuristically we observe that the finite sequence $\left(\frac{\alpha(d,3)}{\overline{\alpha}(d,3)}\right)_{1\le d\le 1000}$ is strictly decreasing. So
the following values might be seen as a good estimate for $\mathcal{P}_3$:

\begin{eqnarray*}
  &&\frac{\alpha(600,3)}{\overline{\alpha}(600,3)}=\frac{2000468396580}{19359502966749}
  \approx 0.103333,\\[1mm]
  &&\frac{\alpha(800,3)}{\overline{\alpha}(800,3)}=\frac{8430487428682}{81665192828999}
  \approx 0.103232,\text{ and }\\[1mm]
  &&\frac{\alpha(1000,3)}{\overline{\alpha}(1000,3)}=\frac{25728695195597}{249377330461249}
  \approx 0.103172 .
\end{eqnarray*}

\section{Orderly generation of integral tetrahedra}
\label{sec_orderly}

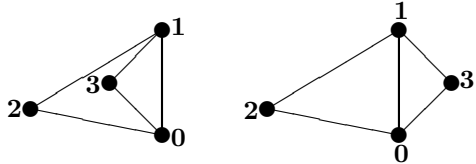
\begin{figure}[ht]
  \begin{center}
    \setlength{\unitlength}{0.7cm}
    \begin{picture}(8.8,2.7)
      \thicklines
      \put(3,0.2){\line(0,1){2}}
      \thinlines
      \put(3,0.2){\line(-1,1){1}}
      \put(2,1.2){\line(1,1){1}}
      \put(3,2.2){\line(-5,-3){2.5}}
      \put(3,0.2){\line(-5,1){2.5}}
      \put(3,0.2){\circle*{0.3}}
      \put(3.17,0){$\mathbf{0}$}
      \put(3,2.2){\circle*{0.3}}
      \put(3.17,2.1){$\mathbf{1}$}
      \put(2,1.2){\circle*{0.3}}
      \put(1.55,1){$\mathbf{3}$}
      \put(0.5,0.7){\circle*{0.3}}
      \put(0.05,0.55){$\mathbf{2}$}
      \thicklines
      \put(7.5,0.2){\line(0,1){2}}
      \thinlines
      \put(7.5,2.2){\line(-5,-3){2.5}}
      \put(7.5,0.2){\line(-5,1){2.5}}
      \put(7.5,0.2){\line(1,1){1}}
      \put(7.5,2.2){\line(1,-1){1}}
      \put(7.5,0.2){\circle*{0.3}}
      \put(7.4,-0.3){$\mathbf{0}$}
      \put(7.5,2.2){\circle*{0.3}}
      \put(7.4,2.4){$\mathbf{1}$}
      \put(8.5,1.2){\circle*{0.3}}
      \put(8.65,1.1){$\mathbf{3}$}
      \put(5,0.7){\circle*{0.3}}
      \put(4.55,0.5){$\mathbf{2}$}
    \end{picture}
  \end{center}
  \caption{Joining two triangles.}
  \label{fig_two_triangles}
\end{figure}

Our strategy to enumerate integral tetrahedra is to merge two triangles along a common side. In Figure \ref{fig_two_triangles} we have depicted the two possibilities in the plane to join two triangles $(0,1,2)$ and $(0,1,3)$ along the side $\overline{01}$. If we 
rotate the triangle $(0,1,3)$ in the $3$-dimensional space from the position on the left in Figure \ref{fig_two_triangles} to the position on the right we obtain tetrahedra and the distance $\delta_{2,3}$ forms an interval $[l,u]$. The restriction to integral 
tetrahedra is fairly easy.

Let us consider the example $\delta_{0,1}=6$, $\delta_{0,2}=\delta_{1,2}=5$, $\delta_{0,3}=4$, and $\delta_{1,3}=3$. Solving $CMD_3=0$ over the positive real numbers yields that the configuration is a 
tetrahedron iff $\delta_{2,3}\in\left(\frac{\sqrt{702-24\sqrt{455}}}{6},\frac{\sqrt{702+24\sqrt{455}}}{6}\right)
\approx(2.297719304,5.806934304)$. Thus there are integral tetrahedra for $\delta_{2,3}\in\{3,4,5\}$. In general we denote such a set of tetrahedra by
$$
  \delta_{0,1},\delta_{0,2},\delta_{1,2},\delta_{0,3},\delta_{1,3},\delta_{2,3}\in[l,r].
$$
This notation permits to implicitly list $\Omega(d^6)$ integral tetrahedra in $O(d^5)$ time.

All integral tetrahedra can be obtained in this manner. So an enumeration method is to loop over all suitable pairs of integral triangles and to combine them. We will go into detail in a while. Before that we have to face the fact that our enumeration method may construct pairs of isomorphic tetrahedra. Looking at Table \ref{table_number_tetrahedra} we see that storing all along the way constructed non-isomorphic integral tetrahedra in a hash table is infeasible. Here we use the concept of orderly generation \cite{winner} which allows us to decide independently for each single constructed discrete structure if we have to take or to reject it. Therefore we have to define a canonical form of an integral tetrahedron. Here we say that a tetrahedron $\mathcal{T}$ with side lengths $\delta_{i,j}$ is canonical if for the lexicographic ordering of vectors $\succeq$,
$$
  (\delta_{0,1},\delta_{0,2},\delta_{1,2},\delta_{0,3},\delta_{1,3},\delta_{2,3})\succeq 
  (\delta_{\tau(0),\tau(1)},\dots,\delta_{\tau(2),\tau(3)})
$$
holds for all permutation $\tau$ of the points $0,1,2,3$. We describe the algorithmic treatment of a canonicity function $\chi(\mathcal{T})\mapsto\{true,false\}$ which decides whether a given integral tetrahedron $\mathcal{T}$ is canonical in Section 
\ref{sec_can}. We have the following obvious lemma:

\begin{lemma}
  If $\chi(\delta_{0,1},\delta_{0,2},\delta_{1,2},\delta_{0,3},\delta_{1,3},\delta_{2,3})=true$ and
  $\chi(\delta_{0,1},\delta_{0,2},\delta_{1,2},\delta_{0,3},\delta_{1,3},\delta_{2,3}+1)=false$ then
  $\chi(\delta_{0,1},\delta_{0,2},\delta_{1,2},\delta_{0,3},\delta_{1,3},\delta_{2,3}+n)=false$ for
  all $n\ge 1$.
\end{lemma}

Thus for given $\delta_{0,1}$, $\delta_{0,2}$, $\delta_{1,2}$, $\delta_{0,3}$, and $\delta_{1,3}$ the possible values for $\delta_{2,3}$ which correspond to a canonical tetrahedron form an interval $[\hat{l},\hat{u}]$. Clearly, the value of $\chi(\delta_{0,1},\delta_{0,2},\delta_{1,2},\delta_{0,3},\delta_{1,3},\delta_{2,3})$ has to be evaluated for $\delta_{2,3}\in\{\delta_{i,j}-1,\delta_{i,j},\delta_{i,j}+1 \mid (i,j)\in\{(0,1),(0,2),(1,2),(0,3),(1,3)\}\}$ only. Thus we can determine the interval $[\hat{l},\hat{u}]$ using
$O(1)$ evaluations of $\chi(\mathcal{T})$.

\begin{algo}
  \label{algo_orderly}
  {\em Orderly generation of integral tetrahedra}\\
  \textit{Input:} Diameter $d$\\
  \textit{Output:} A complete list of canonical integral tetrahedra with diameter $d$\\
  {\bf begin}\\
  \Tsm $\delta_{0,1}=d$\\
  \Tsm {\bf for} $\delta_{0,2}$ {\bf from} $\left\lfloor\frac{d+2}{2}\right\rfloor$ {\bf to} $d$ {\bf do}\\
  \Tsm\Tsm {\bf for} $\delta_{1,2}$ {\bf from} $d+1-\delta_{0,2}$ {\bf to} $\delta_{0,2}$ {\bf do}\\
  \Tsm\Tsm\Tsm {\bf for} $\delta_{0,3}$ {\bf from} $d+1-\delta_{0,2}$ {\bf to} $\delta_{0,2}$ {\bf do}\\
  \Tsm\Tsm\Tsm\Tsm {\bf for} $\delta_{1,3}$ {\bf from} $d+1-\delta_{0,3}$ {\bf to} $\delta_{0,2}$ 
  {\bf do}\\
  \Tsm\Tsm\Tsm\Tsm\Tsm Determine the interval $[\hat{l},\hat{u}]$ for $\delta_{2,3}$\\
  \Tsm\Tsm\Tsm\Tsm\Tsm print $\delta_{0,1},\delta_{0,2},\delta_{1,2},\delta_{0,3},
  \delta_{1,3},[\hat{l},\hat{u}]$\\
  \Tsm\Tsm\Tsm\Tsm {\bf end}\\
  \Tsm\Tsm\Tsm {\bf end}\\
  \Tsm\Tsm {\bf end}\\
  \Tsm {\bf end}\\
  {\bf end}
\end{algo}
We leave it as an exercise for the reader to prove the correctness of Algorithm \ref{algo_orderly} (see \cite{winner} for the necessary and sufficient conditions of an orderly generation algorithm). Since we will see in the next section that we can perform $\chi(\mathcal{T})$ in $O(1)$ the runtime of Algorithm \ref{algo_orderly} is $O(d^4)$. By an obvious modification Algorithm \ref{algo_orderly} returns a complete list of all canonical integral tetrahedra with diameter at most $d$ in $O(d^5)$ time.

We remark that we have implemented Algorithm \ref{algo_orderly} using Algorithm \ref{algo_can_check} for the canonicity check. For the computation of $\alpha(800,3)$ our computer needs only $3.3$ hours which is really fast compared to the nearly $3$ hours needed for a simple loop from $1$ to $\alpha(800,3)$ on the same machine. Due to the complexity of $O(d^4)$ for suitable large $d$ the determination of $\alpha(d,3)$ will need less time than the simple loop from $1$ to $\alpha(d,3)$.

\section{Canonicity check}
\label{sec_can}

In the previous section we have used the canonicity check $\chi(\mathcal{T})$ as a black box. The straight forward approach to implement it as an algorithm is to run over all permutations $\tau\in S_4$ and to check whether $(\delta_{0,1},\dots,\delta_{2,3})\succeq (\delta_{\tau(0),\tau(1)},\dots,\delta_{\tau(2),\tau(3)})$. This clearly leads to running time $O(1)$ but is too slow for our purpose. It may be implemented using $24\cdot 6=144$ integer comparisons. Here we can use the fact that the tetrahedra are generated by Algorithm \ref{algo_orderly}. So if we arrange the comparisons as in Algorithm \ref{algo_can_check} at most $6$ integer comparisons suffice.

\begin{algo}
  \label{algo_can_check}
  {\em Canonicity check for integral tetrahedra generated by Algorithm \ref{algo_orderly}}\\
  \textit{Input:} $\delta_{0,1},\delta_{0,2},\delta_{1,2},\delta_{0,3},\delta_{1,3},\delta_{2,3}$\\
  \textit{Output:} $\chi(\delta_{0,1},\delta_{0,2},\delta_{1,2},\delta_{0,3},\delta_{1,3},\delta_{2,3})$\\
  {\bf begin}\\
  \Tsm{\bf if} $\delta_{0,1}=\delta_{0,2}$ {\bf then}\\
  \Tsm\Tsm{\bf if} $\delta_{0,2}=\delta_{1,2}$ {\bf then}\\
  \Tsm\Tsm\Tsm{\bf if} $\delta_{0,3}<\delta_{1,3}$ {\bf then} {\bf return} $false$\\
  \Tsm\Tsm\Tsm{\bf else}\\
  \Tsm\Tsm\Tsm\Tsm{\bf if} $\delta_{1,3}<\delta_{2,3}$ {\bf then} {\bf return} $false$ {\bf else}
  {\bf return} $true$\\
  \Tsm\Tsm\Tsm{\bf end}\\
  \Tsm\Tsm{\bf else}\\
  \Tsm\Tsm\Tsm{\bf if} $\delta_{1,3}<\delta_{2,3}$ {\bf then} {\bf return} $false$\\
  \Tsm\Tsm\Tsm{\bf else}\\
  \Tsm\Tsm\Tsm\Tsm{\bf if} $\delta_{0,1}=\delta_{0,3}$ {\bf then}\\
  \Tsm\Tsm\Tsm\Tsm\Tsm{\bf if} $\delta_{1,2}<\delta_{1,3}$ {\bf then return} $false$
  {\bf else return} $true$ {\bf end}\\
  \Tsm\Tsm\Tsm\Tsm{\bf else}\\
  \Tsm\Tsm\Tsm\Tsm\Tsm{\bf if} $\delta_{0,1}>\delta_{1,3}$ {\bf then return} $true$\\
  \Tsm\Tsm\Tsm\Tsm\Tsm{\bf else}\\
  \Tsm\Tsm\Tsm\Tsm\Tsm\Tsm{\bf if} $\delta_{1,2}<\delta_{0,3}$ {\bf then return} $false$
  {\bf else return} $true$ {\bf end}\\
  \Tsm\Tsm\Tsm\Tsm\Tsm{\bf end}\\
  \Tsm\Tsm\Tsm\Tsm{\bf end}\\
  \Tsm\Tsm\Tsm{\bf end}\\
  \Tsm\Tsm{\bf end}\\
  \Tsm{\bf else}\\
  \Tsm\Tsm{\bf if} $\delta_{0,2}=\delta_{1,2}$ {\bf then}\\
  \Tsm\Tsm\Tsm{\bf if} $\delta_{0,1}<\delta_{2,3}$ {\bf or} $\delta_{0,3}<\delta_{1,3}$
  {\bf then return} $false$ 
  {\bf else return} $true$ {\bf end}\\ 
  \Tsm\Tsm{\bf else}\\
  \Tsm\Tsm\Tsm{\bf if} $\delta_{0,2}=\delta_{1,3}$ {\bf then}\\
  \Tsm\Tsm\Tsm\Tsm{\bf if} $\delta_{0,3}>\delta_{1,2}$ {\bf or} $\delta_{0,1}<\delta_{2,3}$
  {\bf then return} $false$ 
  {\bf else return} $true$ {\bf end}\\ 
  \Tsm\Tsm\Tsm{\bf else}\\
  \Tsm\Tsm\Tsm\Tsm{\bf if} $\delta_{0,2}=\delta_{0,3}$ {\bf then}\\
  \Tsm\Tsm\Tsm\Tsm\Tsm{\bf if} $\delta_{1,2}<\delta_{1,3}$ {\bf or} $\delta_{0,1}\le\delta_{2,3}$
  {\bf then return} $false$
  {\bf else return} $true$ {\bf end}\\ 
  \Tsm\Tsm\Tsm\Tsm{\bf else}\\
  \Tsm\Tsm\Tsm\Tsm\Tsm{\bf if} $\delta_{0,3}>\delta_{1,2}$ {\bf then}\\
  \Tsm\Tsm\Tsm\Tsm\Tsm\Tsm{\bf if} $\delta_{0,1}\le\delta_{2,3}$ {\bf then return} $true$
  {\bf else return} $false$ {\bf end}\\
  \Tsm\Tsm\Tsm\Tsm\Tsm{\bf else}\\
  \Tsm\Tsm\Tsm\Tsm\Tsm\Tsm{\bf if} $\delta_{0,1}<\delta_{2,3}$ {\bf then return} $true$
  {\bf else return} $false$ {\bf end}\\
  \Tsm\Tsm\Tsm\Tsm\Tsm{\bf end}\\
  \Tsm\Tsm\Tsm\Tsm{\bf end}\\
  \Tsm\Tsm\Tsm{\bf end}\\
  \Tsm\Tsm{\bf end}\\
  \Tsm{\bf end}\\
  {\bf end}
\end{algo}

To prove the correctness of Algorithm \ref{algo_can_check} we consider all vectors $(\delta_{0,1},\dots,\delta_{2,3})$ with $\delta_{i,j}\in\{0,\dots,5\}$, $\delta_{0,1}\ge\delta_{0,2}\ge\delta_{1,2}$, $\delta_{0,2}\ge\delta_{0,3}$, and $\delta_{0,2}\ge\delta_{1,3}$. It suffices to show that Algorithm \ref{algo_can_check} returns the correct value for this finite set of vectors since these inequalities are fulfilled by Algorithm \ref{algo_orderly} and also necessary for $\chi(\mathcal{T})=true$. Algorithm \ref{algo_can_check} can be considered as a binary decision tree. It might be a task to optimize this type of binary decision tree in the worst or in the average case.

\section{Dimensions $\mathbf{m\ge 4}$}
Clearly the question for bounds for $P_m$ arises also for $m\ge 4$. But non-trivial answers seem out of reach by our approach. So far we have no efficient equivalent of Algorithm \ref{algo_can_check} at hand and the number $\alpha(4,d)$ of integral $4$-dimensional simplices with diameter $d$ is $\Omega\left(d^9\right)$. We give the known values of $\alpha(4,d)$ in Table \ref{table_number_simplices}.

\begin{table}[h]
  \begin{center}
    \begin{tabular}{||rr|rr|rr|rr||}
      \hline
      \hline
        $d$ & $\alpha(d,4)$ & $d$ & $\alpha(d,4)$ & $d$ & $\alpha(d,4)$ & $d$ & $\alpha(d,4)$ \\
      \hline  
       1 &       1 & 14 &   12957976 & 27 &   4716186332 & 40 &  162007000505 \\
       2 &       6 & 15 &   24015317 & 28 &   6541418450 & 41 &  202323976907 \\
       3 &      56 & 16 &   42810244 & 29 &   8970194384 & 42 &  251321436143 \\
       4 &     336 & 17 &   73793984 & 30 &  12168243592 & 43 &  310607982160 \\
       5 &    1840 & 18 &  123240964 & 31 &  16344856064 & 44 &  382002253424 \\
       6 &    7925 & 19 &  200260099 & 32 &  21748894367 & 45 &  467627887530 \\
       7 &   29183 & 20 &  317487746 & 33 &  28688094208 & 46 &  569910996879 \\
       8 &   91621 & 21 &  492199068 & 34 &  37529184064 & 47 &  691631229557 \\
       9 &  256546 & 22 &  747720800 & 35 &  48713293955 & 48 &  835911697430 \\
      10 &  648697 & 23 & 1115115145 & 36 &  62769489452 & 49 & 1006370948735 \\
      11 & 1508107 & 24 & 1634875673 & 37 &  80321260053 & 50 & 1207047969441 \\
      12 & 3267671 & 25 & 2360312092 & 38 & 102108730634 & 51 & 1442539675756 \\
      13 & 6679409 & 26 & 3358519981 & 39 & 128999562925 & 52 & 1718015775541 \\
      \hline
      \hline
    \end{tabular}
    \label{table_number_simplices}
    \caption{Number $\alpha(d,4)$ of integral $4$-dimensional simplices with diameter $1\le d\le 52$.}
  \end{center}
\end{table}

\end{document}